\documentclass[a4paper,12pt]{article}

\usepackage[T2A]{fontenc}
\usepackage[utf8]{inputenc}   % older versions of ucs package
\usepackage[russian]{babel}

\usepackage{CJKutf8}

\usepackage{amsmath,amstext,amsfonts,amsthm,amssymb}
\usepackage{eucal}
\usepackage{graphicx}
\usepackage{tikz-cd}
\usetikzlibrary{patterns}
\usetikzlibrary{calc}
\renewcommand{\subsubsection}[1]{\addtocounter{subsubsection}{1}
{\ \\[3pt]\bf \thesubsubsection. \  #1} }
\swapnumbers

{  \theoremstyle{definition}

}
\newcommand{\iso}{\overset{\sim}{\longrightarrow}}
\newcommand{\isom}{\overset{\sim}{=}}
\newcommand{\lra}{\longrightarrow}

\newcommand{\bea}{\begin{eqnarray*}}
\newcommand{\eea}{\end{eqnarray*}}
\newcommand{\bean}{\begin{eqnarray}}
\newcommand{\eean}{\end{eqnarray}}

\newcommand{\bA}{\mathbf{A}}
\newcommand{\bfa}{\mathbf{a}}

\newcommand{\bM}{\mathbf{M}}

\newcommand{\BA}{\mathbb{A}}

\newcommand{\BG}{\mathbb{G}}

\newcommand{\BP}{\mathbb{P}}

\newcommand{\BZ}{\mathbb{Z}}

\newcommand{\nc}{\newcommand}

\nc{\Id}{\text{Id}}
\nc{\la}{\lambda}

\begin{document}

%\centerline{\bf A TORIC DEFINITION OF THE CAYLEY-MOUFANG PLANE}

%\centerline{\bf OCTONIONIC HIRZEBRUCH SURFACES}

\centerline{\bf MOUFANG LOOPS AND TORIC SURFACES}

\

\centerline{\it Picard-Lefschetz reflections in the nonassociative world}

\bigskip\bigskip

\centerline{Vadim Schechtman}

\bigskip\bigskip

\begin{CJK}{UTF8}{min}

%\centerline{ヘッケの縮退代数とコ骸骨\footnote{hecke no shukutai daisu to kogaikotsu}}
%\centerline{ニルヘッケの代数とコ骸骨\footnote{niruhecke no daisu to kogaikotsu}}

\end{CJK}

%\bigskip\bigskip

 \centerline{February 26, 2021}

\

% это продолжение моего сегодняшнего письма; I will keep the notations from it. 

%\begin{tikzpicture}[scale=.4, baseline=(current bounding box.center)]
 
% \node at (0,0){$\bullet$};
%  \draw (0,-4) -- (0,4); 
%  \draw (4,-0) -- (-4,0); 
 % \draw (4,4) --(4, -4) -- (-4,-4) -- (-4,4) -- (4,4); 
%  \node at (-.5, -0.5) {$0$}; 
% \node at (3.5, 3.5){$\hen$};  
%  \end{tikzpicture}
%  \quad\quad 
% { \huge  $\buildrel z^2\over  \lra$}
%   \quad\quad 

\centerline{\bf \S 0. Introduction.}

\

It was establshed in the early 1930-th by the mathematicians of the Hilbert school (Emil Artin, Max Zorn, 
Ruth Moufang) that the notion of a "field" ({\it  K\"orper}) $K$ is the same as the notion 
of a projective plane $\BP^2(K)$: the laws of addition and multiplication in $K$ may be restored 
geometrically from the incidence relations of points and lines in $\BP^2(K)$. 
The associativity of multiplication is equivalent to the Desargue theorem, so that a plane 
where this theorem is true corresponds to
a skew-field ({\it Schiefk\"orper}).

R.Moufang considered a plane with a weakened form of the Desargue theorem - the so called 
{\it der Satz vom volst\"andigen Vierseit $(D_9)$}; the set of nonzero elements 
of the corresponding field ({\it Alternativk\"orper}) carries a nonassociative structure called today a {\it Moufang loop}. 

In this note we propose a construction of some varieties which might be considered 
as "nonsingular toric surfaces" $X_L$\ over a Moufang loop $L$. Suppose that $L = K\setminus \{0\}$ where 
$K$ is a Cayley octonion division algebra over a commutative field $k$. Then, 
as is usual in the toric geometry, 
$X_L$ is patched together out of several copies of affine planes $K^2$. The transition functions 
between various charts have a nice triangular "noncommutative Picard-Lefschetz"\ shape. 
This way we get certain $16$-dimensional varieties over $k$.

We restrict ourselves to some examples. The general case depends on some combinatorial 
"lifting conjecture"\ , see 4.4 below, which might be of interest in itself. 
We note that the gluing construction does not use the full strength of the Moufang identity but 
only a weaker inversion property, see 1.1 below.

I am grateful to B.Toen and J.Tapia for inspiring discussions.

%SAY ABOUT {\it loops with inverse property} SEE [Br], VII, 1.

%SOME NEW $\BP^2$'S ?

\newpage

\

\centerline{\bf \S 1. Moufang loops}

\

{\bf 1.1. Loops.} Recall that a {\it magma} (in the terminology of N.Bourbaki) 
is a set $L$ equipped with an operation
$$
L\times L\lra L,\ (x, y)\mapsto xy
\eqno{(1.1.1)}
$$
which need not to be neither associative nor commutative.

Let us call {\it a loop} a magma $L$ equipped
 with a distingushed element $1$ and a bijection 
$$
L\iso L,\ x\mapsto x^{-1}
%\eqno{(1.2.1)}
$$
such that for all $x, y\in L$

(i) $x1 = x = 1x$;

(ii) 
$$
(x^{-1})^{-1} = 1;
$$ 
$$
x^{-1}x = xx^{-1} = 1;
$$

(iii)
$$
x^{-1}(xy) = (yx)x^{-1} = y,
$$
cf. [CS] 7.1\footnote{where it is called a {\it loop with inversion}.}.

It is shewn in {\it op. cit.} 7.3 that in $L$ 
$$
(yx)y = y(xy)
\eqno{(1.1.2)}
$$
for all $x, y$; we denote this product $yxy$.
It follows that for $a\in \BZ_{\geq 0}$ is defined in $x^a$ unambiguous manner by induction: $x^0 = 1$, 
$$
x^a = x(x^{a-1}) = (x^{a-1})x;
\eqno{(1.1.3)}
$$
we define $x^{-a} := (x^a)^{-1}$. 

{\bf 1.2. Lemma.} {\it In a loop
for all $x, y$
$$ 
(xy)^{-1} = y^{-1}x^{-1}.
\eqno{(1.1.4)}
$$}

{\bf Proof} [H, Lemma 2.11]
$$
(xy)^{-1} = y^{-1}(y(xy)^{-1}) = y^{-1}((x^{-1}(xy))(xy)^{-1}) = y^{-1}x^{-1}. 
$$
$\square$

Note that in this proof the unity $1$ has not been used.

{\bf 1.3. Moufang loops.} Recall that a {\it Moufang loop} is a loop $L$  satisfying 
the axiom:

--- for all $x, y, z$
$$
z(xy)z = zx\cdot zy,
\eqno{(a)}
$$
$$
z(xy) = (zxz)(z^{-1}y),
\eqno{(b)}
$$
$$
(xy)z =(z^{-1})(zyz);
\eqno{(c)}
$$
in fact any of these identities implies two others, 
cf. [CS], 7.4.

%PROVE THAT $(zx)z = z(xz)$

\

{\bf 1.4. Diassociativity.}   
{\bf Artin - Moufang  theorem} (see [Z], Satz 4, p. 127, [M2]) {\it A Moufang loop with two generators is associative.} $\square$

%EXAMPLES: (a) ASSOCIATIVE DIVISION ALGEBRAS

%(b) CAYLEY OCTAVE ALGEBRAS

\

\

\centerline{\bf \S 2. Cayley-Moufang projective plane}

\

{\bf 2.1. Projective line and projective plane in the toric language.} Let $k$ be a field. 

(a) The projective 
line $\BP^1(k)$ with coordinates $(t_0: t_1)$ is covered by two opens both isomorphic to an affine line 
$\BA^1(k) = k$ where we introduce the coordinates
$$
U_0 = \{t_0\neq 0\},\ x := t_1/t_0
$$
and 
$$
U_1 = \{t_1\neq 0\},\ x' := t_0/t_1;
$$
their intersection admits two coordinates
$$
U_{01} = U_0\cap U_1 \isom \BG_m(k) = k^*, x' = x^{-1}.
\eqno{(2.1.1)}
$$
So $\BP^1(k)$ is obtained by gluing of two copies of $\BA^1(k)$ along $\BG_m(k)$, with a transition 
function given by (2.1.1).    

\

(b) The projective 
plane $\BP^2(k)$ with coordinates $(t_0: t_1: t_2)$ is covered by three opens 
$U_i = \{t_i\neq 0\}$ isomorphic to an affine plane $\BA^2(k) = k^2$.

We mark the coordinates in these $3$ charts:

in $U_0$: $x = t_1/t_0, \ y = t_2/t_0$;

in $U_1$: $x' = t_0/t_1, \ y' = t_2/t_1$;

in $U_2$: $x'' = t_0/t_2, \ y'' = t_1/t_2$.

Double intersections $U_{ij} \isom k^*\times k$; transition functions
$$
\phi_{01}(x,y) = (x',y'),\ x' = x^{-1}, \ y' = x^{-1}y
\eqno{(2.1.2a)}
$$
on
$$
U_{01} = \{(x,y)\in U_0|\ x\neq 0\} = \{(x',y')\in U_1|\ x'\neq 0\};
$$
$$
\phi_{12}(x',y') = (x'',y''),\ x'' = y^{\prime -1}x', \ y'' = y^{\prime -1},
\eqno{(2.1.2b)}
$$
on
$$
U_{12} = \{(x',y')\in U_1|\ y'\neq 0\} = \{(x'',y'')\in U_2|\ y''\neq 0\};
$$
$$
\phi_{02}(x,y) = (x'',y''),\ x'' = y^{-1}, \ y'' = y^{-1}x.
\eqno{(2.1.2c)}
$$
on
$$
U_{02} = \{(x,y)\in U_0|\ y\neq 0\} = \{(x'',y'')\in U_2|\ x''\neq 0\}.
$$
On the triple intersection $U_{012} \isom k^{*2}$ we have three coordinate systems:
$(x,y), (x',y')$ and $(x'',y'')$.

{\bf 2.1.1. Cocycle lemma.} {\it 
On the triple intersection $U_{012}$ we have 
$$
\phi_{02} = \phi_{12}\phi_{01}.
\eqno{(2.1.3)}
$$} 

{\bf Proof.} We have
$$
x'' = y^{\prime -1}x' = (y^{-1}x)x^{-1} = y^{-1}
$$
and 
$$
y'' = y'' = y^{\prime -1} = (x^{-1}y)^{-1} = y^{-1}x.
$$
$\square$

Note that in this computation we have only used the identities 1.1 (i) - (iii), no commutativity and 
a weak form of the associativity (1.1.4).

\

{\bf 2.2. Projective plane for a loop.} Let $L$ be a  loop.  We denote {\it an affine line} for $L$  
$$
A = A(L) = L\coprod\{0\}.
$$
We define a {\it projective line} $\BP^1(L)$ as the quotient of a disjoint union
$$
A(L)\times A_\infty(L)
$$
where $A_\infty(L) = L\coprod\{\infty\}$ under an identification
$$
x' = x^{-1}
$$
where $x\in L\subset A(L)$ (resp. $x'\in L\subset A\infty(L)$). Thus
$$
\BP^1(L) = \{0\}\coprod L\coprod \{\infty\}.
$$

We define the {\it projective plane} $\BP^2(L)$ as a quotient of a disjoint union 
of three copies 
$$
A_{0}\coprod A_{1}\coprod A_{2}
$$
where $A_{0}$ (resp. $A_{1}, A_{2}$) is a copy of $A(L)^2 = A\times A$ with coordinates $(x,y)$ 
(resp. $(x',y'), (x'',y'')$) under an equivalence relation:
$$
(x',y') = \phi_{01}(x,y),\ x\neq 0,
$$
$$
(x'',y'') = \phi_{12}(x',y'),\ y'\neq 0,
$$
$$
(x'',y'') = \phi_{02}(x,y),\ y\neq 0.
$$
Here
$$
\phi_{01}:\ L\times A \iso L\times A,
$$
$$
\phi_{12}:\ A\times L \iso A\times L,
$$
$$
\phi_{02}:\ A\times L \iso L\times A.
$$
This set is well defined due to the key 

{\bf 2.2.1. Cocycle relation.} {\it We have
$$
\phi_{02} = \phi_{12}\phi_{01}
$$
on the triple intersection $L\times L$.} 

{\bf Proof} The computation in the proof of 2.1.1 
uses only the identities 1.1 (i) - (iii), no commutativity and 
a weak form of the associativity. $\square$

\

{\bf 2.3. Cellular decomposition.} We have by definition
$$
\BP^2(L) = C_2\coprod C_1\coprod C_0
$$
where
$$
C_2 = U_2 \isom A^2,
$$
$$
C_2 = U_1\setminus U_{10}\isom A,
$$
$$
C_0 = U_2\setminus (U_{20}\cup U_{21})\isom \{*\},
$$
this is what Tits calls a coordinate system on the Cayley projective plane,  
cf. [T], 9.11.2.

\

{\bf 2.4. Examples.} Let $L$ be the variety of nonzero elements of the real Cayley octonions. Then 
$\BP^1_L = S^8$, and $\BP^2_L$ is a $16$-dimensional real variety -  the Cayley projective plane.

\

%EXOTIC (STRANGE) QUATERNIONC PLANES?

\centerline{\bf \S 3. Octonionic toric surfaces}

\

{\bf 3.1. Hirzebruch surfaces.} (a) {\it Commutative version}

 Let $k$ be a field, $a > 0$ be an integer. Following [F], Introduction, p.7, a Hirzebruch surface $F_a$ is  
defined by patching four affine charts $U_i\isom k^2$, $0\leq i\leq 3$ with coordinates and 
transition functions:

$(x,y)$ on $U_0$,

$(x',y') = \phi_{01}(x,y) = (y^{-1}, x)$ on $U_1$,

$(x'',y'') = \phi_{02}(x,y) = (x^{-1}, x^{-a}y^{-1})$ on $U_2$,

$(x''',y''') = \phi_{03}(x,y) = (yx^a, x^{-1})$ on $U_3$.

\

(b) {\it A loop version}

\

Let $L$ be a loop. 
We define $A = A(L)$ as above in \S 2, and consider four copies $U_i\isom A\times A$ with coordinates 
$(x,y), (x',y')$, etc.; inside them six subsets
$$
U_{01} = A\times L\subset U_0,\ U_{10} = L\times A\subset U_1;
$$
$$
U_{12} = A\times L\subset U_1,\ U_{21} = L\times A\subset U_2;
$$
$$
U_{23} = A\times L\subset U_2,\ U_{32} = L\times A\subset U_3;
$$
next
$$
U_{02} = L\times L\subset U_0,\ U_{20} = L\times L\subset U_2;
$$
$$
U_{13} = L\times L\subset U_1,\ U_{31} = L\times L\subset U_3.
$$
next
$$
U_{30} = L\times A\subset U_3,\ U_{03} = A\times L\subset U_0.
$$

Define six patching functions (recall that $x^{\pm a}$ is defined in (1.1.3)):
$$
\phi_{01}: U_{01}\iso U_{10},\ \phi_{01}(x,y) = (y^{-1}, x),
$$
$$
\phi_{12}: U_{12}\iso U_{21},\ \phi_{12}(x',y') = (y^{\prime - 1}, y^{\prime -a}x'),
$$
$$
\phi_{23}: U_{23}\iso U_{32},\ \phi_{23}(x'',y'') = (y^{\prime\prime -1}, x''),
$$
and
$$
\phi_{02}: U_{02}\iso U_{20},\ \phi_{02}(x,y) = (x^{-1}, x^{-a}y^{-1}),
$$
$$
\phi_{13}: U_{13}\iso U_{31},\ \phi_{13}(x',y') = (x^{\prime -1}y^{\prime a}, y^{\prime -1}),
$$
and finally
$$
\phi_{03}: U_{03}\iso U_{30},\ \phi_{03}(x,y) = (yx^{a}, x^{-1}).
$$

{\bf 3.1.1. Claim (cocycle relation).} {\it For all $0\leq i < j < k \leq 3$
$$
\phi_{ik} = \phi_{jk}\phi_{ij}.
$$}

One has to check four relations corresponding to 

$(ijk) = (012), (013), (023)$ and $(123)$. 

This is done directly. 

\

{\bf 3.2. $\BP^2$ blown up at one point.} Classically this variety which we denote $X_1$ may be 
obtained as the toric variety associated with a fan with four generators 
$$
v_0, \ldots, v_3\in \BZ^2
$$
where 
$$
v_0 = (1, 0),\ v_1 = (1, 1),\ v_2 = (0, 1),\ v_3 = (- 1,- 1),  
$$
we number them in the clockwise order,  
cf. [F], 2.4. Thus
$$
X_1 = \cup_{i\in\BZ/4\BZ} U_i
$$
where $U_i$ is the chart $\isom k^2$ corresponding to the dual of the cone with generators $v_i, v_{i+1}$.

Let us describe the corresponding Moufang variety $X_1(L)$ where $L$ is a loop.
It is patched from four charts
$$
U_i(L) \isom A(L)^2
$$
whose coordinates we denote by $(u,v), (u',v')$, etc.

{\bf 3.2.1. Toric formulas.} 
It is convenient to express them in terms of coordinates $x, y$ usual in toric geometry.

The formulas below should be considered as a definition of certain couples of elements 
in the free inversion loop on generators $x, y$; they are the usual toric formulas but the order of factors 
is important.
$$
(u, v) = (x, x^{-1}y),\ 
(u', v') = (y^{-1}x, y)
%\eqno{(3.1.1)}
$$
$$
(u'', v'') = (y^{-1}, y^{-1}x),\ 
(u''', v''') = (x^{-1}y, x^{-1})
\eqno{(3.1.1)}
$$

{\bf 3.2.2.} Returning to our four charts, 
we define inside them the subsets corresponding to double intersections 
$U_{ij}\subset U_i$, $i\neq j$, with 
$$
U_{01}\isom A\times L,\ U_{10}\isom L\times A,
$$
$$ 
U_{12}\isom A\times L, U_{21}\isom L\times A,
$$
$$ 
U_{23}\isom A\times L, U_{32}\isom L\times A,
$$
$$
U_{02} \isom U_{20} = L\times L
$$
$$
U_{13}\isom U_{31}\isom L\times L,
$$
$$
U_{03}\isom L\times A, U_{30}\isom A\times L. 
$$

We define six patching functions:
$$
\phi_{01}:\ U_{01}\iso U_{10},\ \phi_{01}(u,v) = (u',v') = ( v^{-1}, uv);
$$
$$
\phi_{12}:\ U_{12}\iso U_{21},\ \phi_{12}(u',v') = (u'',v'') = (v^{\prime -1},u^{\prime});
$$
$$
\phi_{23}:\ U_{23}\iso U_{32},\ \phi_{23}(u'',v'') = (u''',v''') = 
(v^{\prime\prime -1},v^{\prime\prime -1}u^{\prime\prime});
$$
next
$$
\phi_{02}:\ U_{02}\iso U_{20},\ \phi_{02}(u,v) = (u'',v'') = (v^{-1}u^{-1}, v^{-1});
$$
$$
\phi_{13}:\ U_{13}\iso U_{31},\ \phi_{13}(u',v') = (u''',v''') = (u^{\prime -1},u^{\prime -1}v^{\prime -1});
$$
finally
$$
\phi_{03}:\ U_{03}\iso U_{30},\ \phi_{03}(u,v) = (u''',v''') = (v, u^{-1}).
$$
We need to check four transitivity relations:
$$
\phi_{ik} = \phi_{jk}\phi_{ij}
$$
for $(i, j, k) = (0, 1, 2), (0, 1, 3), (0, 2, 3)$, or  $(1, 2, 3)$. This is done readily. 
One uses only the axiom 1.1 (iii), (1.1.4).

\

{\bf 3.3. Other toric surfaces.} A general nonsingular complete toric surface $X$ is a blowing up at several points either of $F_a$ or of $\BP^2$, cf. [F], 2.5.

Their loop versions are discussed in the next Section.

\

%\centerline{\bf \S 4. Simple (reasonable) cyclic collections}

\centerline{\bf \S 4. Nice cyclic collections}

\

% COMBINATORIAL CONJECTURE

{\bf 4.1. Elementary matrices and pseudoreflections.} Let $a\in\BZ$. We consider four kinds 
of matrices, to be called {\it elementary}, belonging to $SL_2(\BZ)$: 
$$
A(a) = \left(\begin{matrix} 0 & -1\\ 1 & a\end{matrix}\right),\ 
A'(a)  = \left(\begin{matrix} a & 1\\ - 1 & 0\end{matrix}\right),
$$
$$
B(a) = \left(\begin{matrix} -1 & -0\\ a & - 1\end{matrix}\right),\ 
B'(a)  = \left(\begin{matrix} - 1 & a\\ 0 & - 1\end{matrix}\right)
$$
The inverse to an elementary matrix is elementary:
$$
A'(a) = A(a)^{-1},\ 
B(a)^{-1} = B(-a),\ B'(a)^{-1} = B'(-a)
\eqno{(4.1.1)}
$$ 
but the product of elementary matrices 
of different kinds is not elementary.

% WRITE DOWN GOOD PRODUCTS

To each elementary matrix $X(a)$ (for $X\in\{A, A', B, B'\}$) we associate two invertible universal transformations  
$$
X_\epsilon(a):\ L\times L \lra L\times L,\ \epsilon = \pm
\eqno{(4.1.2)}
$$
defined for every loop $L$. 

Namely:
$$
A_+(a)(x,y) = (y^{-1}, xy^a),\ A_-(a)(x,y) = (y^{-1}, y^ax),
$$
$$
A'_-(a)(x,y) = (x^ay, x^{-1}),\ A'_-(a)(x,y) = (yx^a, x^{-1}),
$$ 
$$
B_+(a)(x,y) = (x^{-1}, x^ay^{-1}),\ B_-(a)(x,y) = (x^{-1}, y^{-1}x^a),
$$
$$
B'_-(a)(x,y) = (x^{-1}y^a, y^{-1}),\ B'_-(a)(x,y) = (y^ax^{-1}, y^{-1}).
$$
Their inverses: 
$$
A_\pm(a)^{-1} = A'_{\mp}(a),\ B_\pm(a)^{-1} = B_\mp(-a),\ B'_\pm(a)^{-1} = B'_\mp(-a)
$$
We will say that $X_\epsilon(a)$ is a {\it pseudo-reflection} with matrix $X(a)$, 
or that $X_\epsilon(a)$ is a lifting of $X(a)$. 

We have seen that the inverse to a pseudoreflection is a pseudoreflection.

However the composition of two pseudoreflections is not necessarily a pseudoreflection. 

% GIVE AN EXAMPLE

\ 

{\bf 4.2. Fulton cycles.}
Let $n\geq 3$.  
Let us call a sequence of matrices 
$$
\bA = \{ A(a_i),\ a_i\in\BZ,\  0\leq i\leq n-1\}
%\eqno{(4.2.1)}
$$
a {\it Fulton toric cycle} if 
$$
A(a_0)A(a_1)\ldots A(a_{n-1}) = I.
\eqno{(4.2.1)}
$$
Such a sequence is exactly the data defining a complete nonsingular toric surface, cf. [F], 2.5.

Note that (4.2.1) implies
$$
A(a_i)A(a_{i+1})\ldots A(a_n)A(a_1)\ldots A(a_{i-1}) = I.
%\eqno{(4.2.2)}
$$
for any $0\leq i\leq n-1$. 

Let us identify the set $\{0, \ldots, n-1\}$ with $\BZ/n\BZ$ in the obvious manner which provides us with  
a cyclic order on this set (we imagine its elements as points $p_i$ a circle in the clockwise order).
For each $i, j\in \BZ/nZ$, we define a matrix
$$
A_{ij}(\bfa) = A(a_{j-1})A(a_{j-2})\ldots A(a_i),\ A_{ii}(\bfa) = I
\eqno{(4.2.2)} 
$$
One should imagine $A_{ij}(\bfa)$ as an arrow on the circle
$$
A_{ij}(\bfa): \ p_i \lra p_j
$$
which joins clockwise $p_i$ with $p_j$.

It seems that the following is true.

{\bf 4.2.1. Triangular conjecture.} {\it All matrices $A_{ij}(\bfa)$ are elementary.}

\

{\bf 4.3. A loop version of Fulton cycles.}

{\bf Definition.} {\it A nice cycle of length $n$} is a collection of pseudoreflections 
$$
\bM = \{ M_{ij},\ i, j\in \BZ/n\BZ\},\ 
$$
$$
M_{ij} = X_{\epsilon_{ij}}(a_{ij}),\  X = \{ A, A', B, B'\},
$$ 
such that

(a) For all $i$ $M_{ii} = \Id$, and $M_{i,i+1}$ is of type $A$, i.e. 
$$
M_{i,i+1} = A_{\epsilon_i}(a_i),\ a_i\in \BZ, \epsilon_i = \pm.
$$

(b) For all $i\leq j\leq k$ (in the cyclic order)
$$
M_{ik} = M_{jk}M_{ij}.
$$

Let us denote $M_i := M_{i,i+1}$, and 
$$
\bfa = (a_0, \ldots, a_{n-1})\in \BZ^{\BZ/n\BZ}.
$$

It follows from (b) that
$$
M_0M_1\ldots M_{n-1} = \Id
\eqno{(4.3.1)}
$$
where the order of brackets in the LHS is inessential - the result does not depend on it. More generally
$$
M_{ij} = M_{j-1}M_{j-2}\ldots M_i,
\eqno{(4.3.2)}
$$
for all $i, j$, 
the order of brackets in the RHS being inessential.

The collection of matrices
$$
A(\bfa) = \{A(a_i), \ i\in \BZ/n\BZ\}
$$
is a Fulton cycle.  We will say that $\bM$ is a {\it lifting} of 
$A(\bfa)$.

\

{\bf 4.3. Example: the Hirzebruch cycle ($n = 4$).} We have
$$
A_+(a)A(0)A_-(-a)A(0) = \Id
$$
Here $A(0):= A_+(0) = A_-(0)$. We have:
$$
M_{02} = A_-(-a)A(0) = B_+(-a),\ M_{03} = A'_-(a), \ M_{13} = B'_+(a), 
$$
etc. See 3.1.  

\

{\bf 4.4. Lifting conjecture.} {\it For every Fulton cycle  $A(\bfa)$ there exists a lifting 
$\bM$.}

{\bf 4.5.} 
Given a lifting $\bM$ we can proceed as in \S 3, by gluing  
 the corresponding toric surface $ X(\bM)$ from $n$ copies $U_0, \ldots, U_n$ of $A^2$ and using 
 matrices $M_{ij}$ as patching functions
 $$
 M_{ij}: U_{ij} \iso U_{ji}.
 $$
 
 {\bf 4.6. Mock liftings.} 
 Note that if the initial loop $L$ is a Moufang loop we can proceed as follows  
 to define a lifting. 
 Given a Fulton cycle $A(\bfa)$, 
 we start from transformations $M_i = A_{\epsilon_i}(a_i), 0\leq i\leq n-2$ with arbitrary 
 $\epsilon_i = \pm$. 
 Then define $M_{ij}$ for $0\leq i < j\leq n-2$ by formula (4.3.2): the position of brackets 
 is inessential due to the diassotiativity 1.4. The cocycle condition 4.3 (b) will be satisfied 
 for $0\leq i\leq j\leq k\leq n-2$. These data is enough to glue a variety.
 We can define $M_{n-1}$ by (4.3.1).
 The price we pay for an arbitrary choice of $\epsilon_i$ is that the transformations $M_{ij}$ are in general no longer pseudoreflections. 
 We may call the collection $\bM = \{M_{ij}\}$ a {\it mock lifting} 
 of  $A(\bfa)$. So there are $2^{n-1}$ mock liftings.
   
\bigskip\bigskip

\centerline{\bf References}

\bigskip\bigskip

[A] E.Artin, Coordinates in affine geometry, {\it Rep. Math. Coll. Notre Dame (Indiana)} {\b 2} (1940), 
15-20.

[B] J.Baez, The octonions, {\it Bull. AMS} {\bf 39} (2001), 145-205. 

[Br] R.H.Bruck, A survey of binary systems, {\it Ergebnisse} Bd. {\bf 20}, Springer 1971.

[CS] J.H.Conway, D.A.Smith, On quaternions and octonions: their geometry, arithmetic and symmetry, 
{\it A.K.Peters}, 2003.

[F] W.Foulton, Introduction to toric varieties, {\it Ann. Math. Studies} {\bf 131}, {\it PUP}, 1993.

[H] J.I.Hall, Moufang loops and groups with triality are essentially the same thing, {\it Mem. AMS} 
{\bf 260}, 2019. 

[M1] R.Moufang, Alternativk\"orper and der Satz vom volst\"andigen Vierseit $(D_9)$, {\it Abh. Math. Sem. 
Univ. Hamburg} {\bf 9} (1933), 207-222.

[M2] R.Moufang, Zur Struktur von Alternativk\"orper, {\it Math. Ann.} {\bf 110} (1935), 416 - 430.  

[ST] T.A.Springer, F.D.Veldkampf, Octonions, Jordan algebras and exceptional groups, {\it Springer}, 2000.

[T] J.Tits, Buildings of spherical type and finite BN-pairs, {\it Lecture Notes in Math.} {\bf 386}, 
{\it Springer}, 1974.

[Z] M.Zorn, Theorie der alternativen Ringe, {\it Abh. math. Seminar Univ. Hamburg} {\bf 8} (1931), 
123-147.

%  WRITE HERE

\end{document}